\documentclass{amsart}
\usepackage{amssymb,amsthm,amscd}

\newcommand{\F}{{\mathbb F}}
\newcommand{\N}{{\mathbb N}}
\newcommand{\Q}{{\mathbb Q}}
\newcommand{\C}{{\mathbb C}}
\newcommand{\Z}{{\mathbb Z}}
\newcommand{\fa}{{\mathfrak a}}
\newcommand{\fm}{{\mathfrak m}}
\newcommand{\fp}{{\mathfrak p}}
\newcommand{\cC}{{\mathcal C}}
\newcommand{\cO}{{\mathcal O}}
\newcommand{\lra}{\longrightarrow}
\newcommand{\too}{\longmapsto}
\newcommand{\la}{\langle}
\newcommand{\ra}{\rangle}
\newcommand{\eps}{\varepsilon}
\newcommand{\Cl}{\operatorname{Cl}}

\title{Jacobi and Kummer's Ideal Numbers}
\author{Franz Lemmermeyer}
\address{M\"orikeweg 1, 73489 Jagstzell}
\email{hb3@ix.urz.uni-heidelberg.de}

\begin{document}

\begin{abstract}
In this article we give a modern interpretation of Kummer's
ideal numbers and show how they developed from Jacobi's
work on cyclotomy, in particular the methods for studying
``Jacobi sums'' which he presented in his lectures on number 
theory and cyclotomy in the winter semester 1836/37.
\end{abstract}

\maketitle

\begin{center}
{\em Dem Andenken an Herbert Pieper\footnote{Eine W\"urdigung von 
E. Knobloch zu Piepers 65. Geburtstag findet sich auf \\
{\tt http://www.uni-potsdam.de/u/romanistik/humboldt/hin/hin16/knobloch.htm}.}} 
(1943 -- 2008) {\em gewidmet.}
\end{center}

\bigskip

Every mathematician nowadays is familiar with the notion of an 
ideal in a ring. Ideals were introduced by Dedekind when he 
generalized Kummer's ideal numbers to general number fields. 
Kummer had invented ideal numbers in order to restore some 
kind of unique factorization in cyclotomic rings, and in the 
literature one usually finds the following characterizations
of Kummer's invention:
\begin{enumerate}
\item Kummer's idea was brilliant and new; there were no traces of
      it in the number theoretical work of his predecessors: it
      appeared out of the blue and solved the ``problem'' of
      nonunique factorization in a way reminiscent of Alexander
      the Great's solution of the Gordian knot.
\item Kummer's definition of an ideal prime is difficult to understand
      and not easy to use in practice. 
\end{enumerate}
These opinions seem to be generally accepted: Dickson \cite{Dick}, 
to quote a typical example, wrote in connection with his review of 
Reid's textbook on algebraic number theory:
\begin{quote}
He [Reid] wisely did not attempt to give any idea of Kummer's
ideal numbers, the operations on which are so delicate that
one must use the utmost circumspection (as remarked by Dedekind
in his important historical papers in Darboux's Bulletin).
\end{quote}
This explains why there are hardly any expositions of Kummer's 
theory of ideal numbers; among the few exceptions are Edwards'
book \cite{EdwF} and Soublin's article \cite{Sou}. Nevertheless, 
Kummer's articles on the arithmetic of cyclotomic number fields 
are fairly easy to read\footnote{The majority of Kummer's articles 
is written in German, only a few in Latin and French.} by simply 
replacing the expression ``ideal number'' by the word ``ideal''. 
This works very well except for his work on the foundations of 
his theory of ideal numbers: the problems Kummer had to overcome 
when he introduced ideal numbers seem to differ fundamentally 
from the obstacles Dedekind had to deal with when he introduced 
ideals.

Let us now summarize the content of this article:
\begin{itemize}
\item We start by giving a brief summary of Jacobi's
      lectures on number theory and cyclotomy in the winter
      semester 1836/37; as we will show, they played a major
      role in the development of Kummer's notion of ideal numbers.
\item Next we address the role of Fermat's Last Theorem in Kummer's
      early work.
\item Afterwards we explain the ``Jacobi maps'', certain substitutions
      used by Jacobi that turned out to be the key idea later used
      by Kummer when he invented ideal numbers. 
\item In our discussion of Kummer's ideal numbers we offer a translation 
      of ideal numbers into the modern mathematical language which is 
      more faithful than the simple substitution of ``ideal'' for  
      ``ideal number''. We try to correct the historical picture of 
      the development of Kummer's ideal numbers\footnote{See B\"olling 
      \cite{Boe}, Edwards \cite{EdwF,Edw1,Edw2,Edw3,Edw4}, Neumann 
      \cite{Neu} and Soublin \cite{Sou}.} by showing that the notion 
      of ideal numbers used by Kummer is perfectly natural, and that it 
      is based to a large degree on ideas put forth by Jacobi in his 
      investigations in cyclotomy. Moreover, a theory of divisibility 
      built on these ideas is hardly more complicated than Dedekind's 
      approach; Jung, in his introductory lectures \cite{Jung} on the 
      arithmetic of quadratic number fields, uses an approach that is 
      very close to Kummer's first attempt at defining ideal numbers, 
      and so do Stevenhagen \& Lenstra in \cite[p. 33]{SL}.
\item Finally we discuss the relevance of the notion of integral closure
      for Kummer's work by looking carefully at the concept of singularity 
      in number theory and algebraic geometry. 
\end{itemize} 
It seems that the importance of Jacobi's K\"onigsberg lectures \cite{Jac}
on number theory for the development of ideal numbers has not been 
recognized before. Apparently Kummer had carefully studied a copy of 
Jacobi's lectures. We know from Jacobi himself that Kummer had access
to these lecture notes: In 1846, Jacobi had a note from 1837 reprinted 
in Crelle's Journal \cite{JacKrT} and added a footnote in which he said:
\begin{quote}
Diese aus vielfach verbreiteten Nachschriften der oben
er\-w\"ahn\-ten Vorlesung (an der K\"onigsberger Universit\"at) 
auch den Herren Dirichlet und Kummer seit mehreren Jahren 
bekannten Beweise sind neuerdings von Hrn. Dr. Eisenstein im 
27ten Bande des Crelleschen Journals auf S. 53 publicirt 
worden.\footnote{These proofs, which were known for several 
years to Mr. Dirichlet and Mr. Kummer, among others, through 
widely circulated notes of the lectures (at the University of 
K\"onigsberg) mentioned above, have lately been published by 
Dr. Eisenstein on p. 53 of the 27th volume of Crelle's Journal.}
\end{quote}
Whether Eisenstein really had seen Jacobi's lecture notes prior to
1846 is open to debate; most (if not all) historians of mathematics 
seem to agree that Eisenstein developed his proofs of quadratic, 
cubic and quartic reciprocity laws independently from Jacobi. 
Later, both Kummer and Eisenstein employed (without attribution)
some form of $p$-adic development of the logarithm, which first 
appeared in Jacobi's lectures.

\section{Jacobi's K\"onigsberg Lectures}

In his lectures \cite{Jac} on number theory during the winter 
semester 1836/37, Jacobi introduced his audience to the basic 
theory of Gauss and Jacobi sums (at that time, this was called 
the theory of cyclotomy (Kreisteilung)), and applied these techniques 
to derive the quadratic, cubic and quartic reciprocity laws,  as well 
as results on the representation of primes by quadratic forms that
led him to conjecture Dirichlet's class number formula for 
binary quadratic forms of negative prime discriminants. 

Today, these results are proved using the ideal decomposition of 
Gauss and Jacobi sums. The fact that Jacobi's proofs are 
essentially equivalent to the modern proofs implies that he must
have possessed a technique that allowed him to express the 
essential content of the prime ideal factorizations of Gauss and
Jacobi sums. We will say more on this in Sections \ref{SJM} and
\ref{SMI}, and now turn to the content of Jacobi's lectures.

In the first 5 lectures, Jacobi presents elementary number theory:
congruences, primes, Euler's phi function, the theorem of Euler-Fermat, 
quadratic residues, and the Legendre symbol. The second part deals 
with cyclotomy: he introduces Gauss and Jacobi sums (without their 
modern names, of course), develops their basic properties, and 
explains the connections between Gauss sums and certain binomial 
coefficients; in particular, he proves Gauss's famous congruence 
$2a \equiv \pm \binom{2n}{n} \bmod p$ for primes $p = a^2 + 4b^2 = 4n+1$. 
The third part of the lectures is dedicated to applications of cyclotomy 
to number theory. Jacobi begins with his proof of the quadratic 
reciprocity law (the one Legendre included in his Th\'eorie des Nombres,
and which led to the priority dispute between Jacobi and Eisenstein
that began with Jacobi's footnote quoted above), then derives results 
by Dirichlet on quartic residues as well as the full quartic 
reciprocity law. After briefly discussing cubic residues and 
proving the cubic reciprocity law, Jacobi then shows that if 
$\lambda \equiv 3 \bmod 4$ and $p = \lambda n + 1$ are primes (this 
implies that $(\frac{-\lambda}p) = +1$), then $4p^h = x^2 + \lambda y^2$
for integers $x, y$ and some positive integer $h$ that can be expressed
as a sum of Legendre symbols, and which Jacobi conjectured to be
equal to the number of classes of binary quadratic forms with 
discriminant $-\lambda$. In the last few lectures, he deals with 
similar problems for primes $\lambda \equiv 1 \bmod 4$. 

It is remarkable how quickly Jacobi led his audience from the basic 
facts of elementary number theory right into the middle of research 
problems he was working on.\footnote{Many of the results that Jacobi 
presented in his lectures had also been obtained around 1830 by Cauchy, 
who published his theory in a long series of articles \cite{Cauchy} in 
1840. Lebesgue \cite{Leb} later gave simplified proofs for the main 
results of Jacobi and Cauchy ``on just a few pages'', as he proudly 
remarked.}

Jacobi's principal technique for studying Gauss and Jacobi sums were 
certain substitutions, whose key role in the era before Kummer's ideal 
numbers was emphasized by Frobenius \cite[p. 117--118]{Rei}:
\begin{quote}
Als Cauchy, Jacobi und Kummer angefangen hatten, die Untersuchungen
von Gauss \"uber complexe Zahlen auf allgemeinere aus Einheitswurzeln
gebildete algebraische Zahlen auszudehnen, ergab sich das unerw\"unschte
Resultat, da\ss{} in diesem Gebiete zwei Zahlen nicht immer einen
gr\"o\ss{}ten gemeinsamen Divisor besitzen, und da\ss{} Producte 
unzerlegbarer Factoren einander gleich sein k\"onnen, ohne da\ss{}
die Factoren einzeln \"ubereinstimmen. Die Gleichheit solcher Producte
konnte man daher immer nur durch besondere Kunstgriffe beweisen,
zu denen namentlich der geh\"orte, durch Substitution gewisser
rationalen Zahlen f\"ur die algebraischen die untersuchten Gleichungen
in Congruenzen zu verwandeln. Mit den Methoden, solcher Schwierigkeiten
Herr zu werden, besch\"aftigt sich auch ein gro\ss{}er Theil von
Kronecker's Dissertation.\footnote{When Cauchy, Jacobi and Kummer 
started to extend Gauss's investigations on complex numbers to 
general algebraic numbers formed with roots of unity, they came
across the unpleasant fact that in this domain two numbers do not 
always have a greatest common divisor, and that products of irreducible
factors can be equal without the factors being the same. The equality
of such products could be proved only by resorting to certain tricks,
notably the one that turns the equations under investigation into
congruences by substitutions of certain rational numbers for the 
algebraic numbers. A large part of Kronecker's dissertation deals 
with the methods for mastering such problems.}
\end{quote}
We remark that the first part of Kronecker's thesis (written under
the supervision of Kummer) deals with the basic arithmetic of cyclotomic 
fields. The following result is related to the finiteness of 
the class number of cyclotomic number fields once Kummer had introduced
ideal numbers and the class group: Let $M_\lambda$ denote the maximum 
of the norm of $x + x_1\eps + \ldots + x_{\lambda-1} \eps^{\lambda-1}$
as $-1 \le x_j \le 1$ (here $\eps$ is a primitive $\lambda$th root
of unity). Then for any prime $p$ there is an integer $n < M_\lambda$ 
such that $np$ is a norm from $\Q(\eps)$. Kronecker observes that this
is analogous to the finiteness of the number of reduced forms. In the
second part of his thesis, Kronecker proved ``Dirichlet's'' unit theorem 
for cyclotomic extensions.

\section{Kummer and Fermat's Last Theorem}

The story according to which Kummer, at the beginning of his
career, gave a proof of Fermat's Last Theorem in which he
erroneously assumed unique factorization in the number rings
$\Z[\alpha]$ of $\lambda$-th roots of unity\footnote{In Jacobi's
and Kummer's notation, $\lambda$ is an odd prime number and $\alpha$
a primitive root of the equation $\alpha^\lambda = 1$. The ring
$\Z[\alpha]$ consists of all $\Z$-linear combinations of powers
of $\alpha$.} probably first appeared in the ``Ged\"achtnisrede'' 
on Kummer by Hensel \cite{Rei}, and is now believed to be false. 
Indeed, Edwards \cite{EdwF,Edw1,Edw2}, Neumann \cite{Neu} and 
B\"olling \cite{Boe2} have shown that Kummer's first article on 
cyclotomy dealt with the factorization of primes $\lambda m +1$ 
in the rings $\Z[\alpha]$ of $\lambda$-th roots of unity, and 
that his (false) result implied unique factorization in $\Z[\alpha]$. 

Hensel \cite[p. 93]{Rei} even claimed that Kummer developed his 
theory of ideal numbers only because of Fermat's Last Theorem, and 
this is definitely not true (see also \cite{Corry}). The weaker 
claim that Kummer invented his ideal numbers in connection with 
his work on Fermat's Last Theorem is a story that perhaps originated 
in a short note by Kronecker \cite{Kro}, who claimed
\begin{quote}
So f\"uhrte das Reciprocit\"atsgesetz f\"ur quadratische Reste 
schon zur weiteren Ausbildung der Theorie der Kreistheilung, und 
der be\-r\"uhm\-te Fermat'sche Satz gab Hrn. Kummer vor etwa dreissig 
Jahren die haupts\"achlichste Anregung zu jenen von so gl\"ucklichem
Erfolge gekr\"onten Untersuchungen, auf denen das Reuschle'sche Werk
basirt und deren Weiterf\"orderung es zugleich gewidmet 
ist.\footnote{Thus the reciprocity law for quadratic residues led 
to the further development of the theory of cyclotomy, and Fermat's
famous theorem was thirty years ago Kummer's main motivation for his 
successful investigations on which the work of Reuschle is built and 
to whose further advancement it is dedicated.}
\end{quote}

Kummer himself left no doubt as to which problems motivated his 
work; in \cite{Kum}, he writes
\begin{quote}
Es ist mir gelungen, die Theorie derjenigen complexen Zahlen, 
welche aus h\"oheren Wurzeln der Einheit gebildet sind und welche
bekanntlich in der Kreistheilung, in der Lehre von den Potenzresten
und den Formen h\"oherer Grade eine wichtige Rolle spielen, zu
vervollst\"andigen und zu vereinfachen; und zwar durch die Einf\"uhrung
einer eigenth\"umlichen Art imagin\"arer Divisoren, welche ich 
{\em ideale complexe Zahlen} nenne\footnote{I have succeeded in 
completing and simplifying the theory of those complex numbers
formed from higher roots of unity, which play an important role
in cyclotomy, the theory of power residues, and of the forms of higher
degrees; this was accomplished through the introduction of a special
kind of imaginary divisors, which I call {\em ideal complex numbers};};
\end{quote}

Thus in 1845, right after he had worked out the basic theory of
ideal numbers, Kummer mentions the theory of cyclotomy and power
residues as the driving force behind his work. But when Kronecker, 
who had turned from Kummer's pupil to his closest friend, says
that applications to Fermat's Last Theorem had been on Kummer's
mind during his work on ideal numbers, we cannot simply dismiss 
this as nonsense. In fact, Kummer told Kronecker, in his letter from 
April 2, 1847, that he has found a proof of Fermat's Last Theorem for 
a certain class of exponents, and writes:
\begin{quote}
Der obige Beweis ist erst drei Tage alt, denn erst nach Beendigung 
der Recension\footnote{Kummer is talking about his review of the 
first volume of Jacobi's collected works.} fiel es mir ein wieder 
einmal diese alte Gleichung vorzunehmen, und ich kam diesmal bald 
auf den richtigen Weg.\footnote{This proof is only three days old, 
since only after finishing the review I took up this old
equation again, and this time I soon found the right approach.}
\end{quote}
We know from Kummer's letters that he had discussed his ideal numbers
with both Dirichlet and Jacobi, and although there cannot be any
doubt that Kummer drew his motivation for introducing ideal numbers 
from Jacobi's work on Jacobi sums, it is hard to believe that Dirichlet, 
who had proved Fermat's Last Theorem for the exponents $n = 5$ and 
$n = 14$, failed to point out possible applications of unique 
factorization into ideal primes to the solution of Fermat's Last Theorem. 
The fact that Kummer never mentioned Fermat's Last Theorem as a possible 
motivation for his theory of ideal numbers before 1847 is also not 
surprising; if he was convinced that such an application was possible, 
it was natural for him to keep it to himself until he had worked out 
a proof.

Kummer found such a proof (for primes satisfying certain conditions,
which, as he later showed, hold for regular primes) at the end of 
March 1847. At the beginning of March 1847, Dirichlet received a
letter from his friend Liouville concerning Lam\'e's attempted proof 
of Fermat's Last Theorem; Liouville inquired what Dirichlet knew
about unique factorization in cyclotomic rings of integers. It is 
not known whether Kummer had seen this letter before he took up the 
Fermat equation at the end of March (I think he did); what we do know 
is that Kummer sent his proof to Dirichlet on April 11, and that he
sent a letter to Liouville on April 28. Kummer had already written 
a short note \cite{KumFLT} on the Fermat equation 
$x^{2p} + y^{2p} = z^{2p}$ in 1837, but he had only used elementary 
means there. Kronecker's statements make it plausible that Kummer 
was well aware of possible applications of the theory of cyclotomy to 
the proof of Fermat's Last Theorem, and that he had looked at this 
problem occasionally while he developed his theory of ideal numbers. 
Nevertheless, Fermat's Last Theorem did apparently not play a 
decisive role in Kummer's work before Liouville's letter from March
1847.

Kronecker also knew about Fermat's Last Theorem very early on:
in fact he had used the following claim as the last thesis for 
his disputation at his Ph.D. defense: 
\begin{quote}
Fermatius theorema suum inclytum non demonstravit.\footnote{Fermat 
did not prove his famous theorem.}
\end{quote}

\section{Quadratic Forms vs. Quadratic Number Fields}

In this section we will address two puzzling questions that
do not seem to be related at first:
\begin{enumerate}
\item Both Kummer and Kronecker knew that Kummer's theory of ideal 
      numbers, when applied to numbers of the form $a+b\sqrt{m}$, 
      is essentially equivalent to Gauss's theory of quadratic 
      forms of discriminant $4m$. Yet it was Dedekind who worked out 
      this theory, and he did so as late as 1871! Why wasn't this
      done sooner? 
\item It is clear from the reactions of Jacobi, Dirichlet and 
      Eisenstein to Kummer's retraction of his manuscript that all
      three were fully aware of the failure of unique factorization 
      in cyclotomic (and probably also in quadratic) rings of integers. 
      Why were they all silent on this topic?
\end{enumerate} 
It was suggested that the reason why e.g. Kummer did not develop
a theory of ideal numbers in quadratic fields was the problem
coming from Kummer's choice of the ring $\Z[\sqrt{m}\,]$; this 
is not the maximal order in $\Q(\sqrt{m}\,)$, and we will see
below that this makes his theory of ideal primes break down. This 
argument, however, is not fully convincing: in his proof of the 
$p$-th power reciprocity law, Kummer had to study ideal classes 
in certain orders of Kummer extensions $\Q(\zeta_p,\sqrt[p]{\mu}\,)$,
and these orders were also not maximal. Kummer avoided the problems
caused by the primes dividing $p\mu$ by simply excluding them and
restricting his attention to elements coprime to $p\mu$. 

Dirichlet, in his proofs of Fermat's Last Theorem for the exponents 
$n=5$ in 1828 \cite{Dir0} and $n=14$ in 1832 \cite{Dir1}, did 
use algebraic numbers of the form $a+b\sqrt{5}$ and $a+b\sqrt{-7}$ 
for integers $a, b$, but for deriving their basic properties he 
employed the theory of quadratic forms. Also in 1832, Gauss 
published his second memoir \cite{Gau2} on the theory of 
biquadratic residues. In this article, Gauss proved that 
$\Z[i]$ is factorial and that $\Z[i]$ is Euclidean, but the 
proof of unique factorization is not based on the Euclidean 
algorithm but rather on the fact that the binary quadratic forms with
discriminant $-4$ have class number $1$. Dirichlet, in his article 
\cite{Dir2} on the quadratic reciprocity law in $\Z[i]$, essentially 
copies Gauss's proof. Only ten years later Dirichlet remarked
that domains with a Euclidean algorithm have unique factorization.

The reason why Gauss, Dirichlet, Jacobi and Eisenstein preferred 
the theory of forms was that this was a perfectly general 
theory, whereas arguments based on unique factorization only worked 
in very special cases.  This also seems to be the primary 
reason why the theory of ideal numbers in quadratic number fields 
was not seen as an important contribution to mathematics; the 
raison d'\^etre of ideal numbers was their role in proofs of 
reciprocity laws and Fermat's Last Theorem, which were based on
cyclotomic number fields, and quadratic number fields were not yet 
studied for their own sake. Dirichlet's class number formula was 
stated in terms of binary quadratic forms, and even the class 
number formula for quadratic extensions of $\Q(i)$ was proved 
using the language of quadratic forms with coefficients from $\Z[i]$. 
But when Dirichlet learned from Kummer that he had found a substitute 
for unique factorization in general cyclotomic fields, he must have 
realized the potential of this theory right away.

\section{Jacobi Maps}\label{SJM}
After Gauss had given two proofs (his fourth and sixth) of the 
quadratic reciprocity law using Gauss sums, it became clear that 
their generalization was the key to higher reciprocity 
laws\footnote{Eventually, however, it turned out that the
decomposition of Gauss sums only gives a piece of the reciprocity
law for $\ell$-th powers, namely Eisenstein's reciprocity law.
This is enough to derive the full version for cubic and quartic
residues, but for higher powers, Kummer had to generalize Gauss's
genus theory from quadratic forms to class groups in Kummer 
extensions of cyclotomic number fields.}. 
These Gauss and Jacobi sums for characters of higher order were 
studied by Jacobi\footnote{Cauchy also studied these sums, but 
his lack of understanding higher reciprocity kept him from going 
as far as Jacobi did.}, who first collected their basic properties. 

In order to describe Jacobi's results on Jacobi sums, let us first 
explain his notation. Let $p = m\lambda+1$ denote a prime number, 
$g$ a primitive root modulo $p$, $\alpha$ a root of the equation 
$\alpha^\lambda = 1$, and $x$ a root of the equation $x^p = 1$. Then 
$$ (\alpha,x) = x + \alpha x^g + \alpha^2 x^{g^2} + \ldots
                  + \alpha^{p-2} x^{g^{p-2}} $$
is a ``Gauss sum'', whose $\lambda$-th power does not depend on $x$:
$$ (\alpha,x)^\lambda \in \Z[\alpha]. $$
When $x$ is fixed, or when he is studying expressions like 
$(\alpha,x)^\lambda$ that do not depend on $x$, Jacobi 
simply writes $(\alpha) = (\alpha,x)$.

Now assume that $r$ is a primitive root of the equation 
$r^{p-1} = 1$ (put $\alpha = r$ and $m=1$ above). Then 
$$ \psi r = \frac{(r^i)(r^k)}{(r^{i+k})} $$
is a ``Jacobi sum'', which also is independent of $x$. Its main
properties are
$$ \psi r \cdot \psi(r^{-1}) = p, $$
as Jacobi proved in the XIIIth lecture, and the fact that 
$\psi r$ is an element of $\Z[r]$ (see \cite[XXVI. Vorl.]{Jac}):
\begin{quote}
Die Funktion $\psi r$ besteht aus ganzen positiven Zahlen, welche in 
die verschiedenen Potenzen von $r$ multiplicirt sind,\footnote{The
function $\psi r$ consists of positive integers, multiplied by
the different powers of $r$,}
\end{quote}
Jacobi then continues:
\begin{quote}
Wir wollen hier f\"ur $r$ eine ganze Zahl $g$ setzen, welche
primitive Wurzel der Kongruenz $g^{p-1} \equiv 1 \pmod p$ ist.
Dadurch \"andern sich unsere Gleichungen nur so, dass sie 
Kongruenzen in bezug auf den Modul $p$ werden.\footnote{We now 
want to substitute an integer $g$ for $r$, which is a primitive
root of the congruence $g^{p-1} \equiv 1 \pmod p$. This will
change our equations only in so far as they now become congruences
with respect to the modulus $p$.}
\end{quote}
That ``equations become congruences'' is Jacobi's way of expressing 
the fact that this substitution commutes with addition and 
multiplication, so we can translate it into modern language as 
\begin{quote}
``the substitution $r \too g \bmod p$ induces a ring homomorphism 
  $\phi_g:\Z[r] \lra \Z/p\Z$''.
\end{quote}

Jacobi now writes $\psi g$ for the residue class mod $p$ he gets
by substituting $g \bmod p$ for $r$, and proves the fundamental 
congruence
\begin{equation}\label{FC}
 \psi g \equiv \begin{cases}
        \quad 0   & \pmod p \quad \text{ if } i + k < p-1, \\
      \frac{(2(p-1)-i-k)!}{(p-1-i)! (p-1-k)!} 
                  & \pmod p \quad \text{ if } i + k > p-1,
    \end{cases} 
\end{equation}
where $i, k$ denote integers $0 < i, k < p-1$.

For studying e.g. Jacobi sums of order $\lambda$ for some
$\lambda \mid (p-1)$, Jacobi writes $p-1 = m\lambda$ and replaces
$r$ by $\alpha = r^m$ and $g$ by $g^m$. Jacobi's fundamental 
congruence (\ref{FC}) then shows that $\psi g^m$ is divisible
by $p$ if and only if $im+km < \lambda$.

\section{A Modern Interpretation of Jacobi's results}\label{SMI}

Let us now see why Jacobi's fundamental congruence (\ref{FC}) 
implies the prime ideal factorization of Jacobi sums (in the 
following, we assume familiarity with basic properties of
character sums; see e.g. \cite{BEW,Lem}).

Consider a prime $p = m\lambda+1$, a primitive root $g \bmod p$,
and a primitive $\lambda$-th root of unity $\zeta$. Then there is a unique
character $\chi$ of order $\lambda$ (a surjective homomorphism of
$(\Z/p\Z)^\times \lra \la \zeta \ra$) such that $\chi(g) = \zeta$.
Jacobi's functions
$$\psi \alpha = \psi_{i,k}\alpha 
              = \frac{(\alpha^i)(\alpha^k)}{(\alpha^{i+k})}$$
for $\lambda$-th roots of unity $\alpha = r^m$ then correspond to 
our Jacobi sums 
$$ J(\chi^i,\chi^k) = - \sum_{t=1}^{p-1} \chi^i(t) \chi^k(1-t)
                    = - \psi_{i,k} \alpha $$ 
(this is the sign convention used in \cite{Lem}; in the notation
used by \cite{BEW}, we have $J(\chi^i,\chi^k) = \psi_{i,k}(\alpha$))
with $J(\chi^i,\chi^k) \in \Z[\zeta]$. 

As a special case consider the character $\chi$ of order $\lambda = 4$ 
defined by the quartic residue symbol $\chi = [\frac{\cdot}{\pi}]$
modulo the prime $\pi = a+bi$ with norm $p$ in $\Z[i]$. We have 
$\chi(g) \equiv g^{(p-1)/4} \bmod \pi$, and we can choose the
primitive root $g \bmod p$ in such a way that $\chi(g) = i$.
With these normalizations, we find 
$\psi \alpha = -J(\chi,\chi)$ for $i = k = \frac{p-1}4$, and 
$\psi \alpha = -J(\chi^3,\chi^3)$ for $i = k = \frac{3(p-1)}4$.
Jacobi's congruence implies that $J(\chi,\chi) \equiv 0 \bmod \pi$
and $\pi \nmid J(\chi^3,\chi^3)$ (the actual congruence gives,
as Jacobi observes, a proof of Gauss's result that, for primes 
$p = 4m+1$, we have $p = a^2 + b^2$, where the odd integer $a$ is 
determined (up to sign) by the congruence 
$a \equiv \frac12 \binom{2m}{m} \bmod p$).
Since Jacobi sums have absolute value $\sqrt{p}$, this implies the 
``prime ideal factorization'' $(J(\chi,\chi)) = (\pi)$ of the 
quartic Jacobi sum.

In general, write $p = m \lambda +1$, let $\zeta = r^m$ be a 
primitive $\lambda$-th root of unity, and consider a character 
of order $\lambda$ on $\Z[\zeta]$ with $\chi(g) = \zeta$. Letting 
$\fp$ denote the prime ideal $\fp = (p, \zeta - g^m)$ we find 
$(\frac g\fp) \equiv g^m \equiv \zeta \bmod \fp$ and thus
$(\frac g\fp) = \zeta$ for the $\lambda$-th power residue symbol 
$(\frac{\cdot}{\fp})$. With this choice of $\fp$, we therefore 
have $\chi = (\frac{\cdot}{\fp})$.

Jacobi's congruence (\ref{FC}) then shows $\fp \mid J(\chi^t,\chi^t)$ 
if and only if $0 < 2t < \lambda$. Let $\sigma_j$ denote the 
automorphism $\zeta \too \zeta^j$ of $\Q(\zeta)/\Q$; then 
$\fp \mid J(\chi^t,\chi^t)$ if and only if 
$\sigma_t^{-1}\fp \mid J(\chi,\chi)$. This implies that 
$J(\chi,\chi)$ is divisible by all the prime ideals
$\sigma_t^{-1}\fp$ with $0 < 2t < \lambda$. Since these prime
ideals are pairwise disjoint, and since $J(\chi,\chi)$ is
an algebraic integer with absolute value $\sqrt{p}$, we conclude 
that 
$$ (J(\chi,\chi)) = \fp^s, \quad 
      s = \sum_{t=1}^{\lambda-1} 
          \Big\lfloor \frac{\lambda-2t}{\lambda}\Big\rfloor \sigma_t^{-1} $$
is the complete prime ideal factorization of the Jacobi sum
$J(\chi,\chi)$. Replacing $t$ by $-t$ in this summation gives
\cite[Cor. 11.5]{Lem}.

Jacobi maps (replacing roots of unity in $\C$ with roots of unity
in $\F_p$) were a tool that allowed Jacobi to state and
prove results that we would describe using prime ideals. As we will
see in the next section, Jacobi maps were indeed used by Kummer for
his first tentative definition of ideal numbers. Not only that, 
after Kummer had worked out the factorization of these Jacobi sums 
into ideal prime factors, he even remarked \cite[p. 362]{KumZ} 
\begin{quote}
Dieses Jacobische Resultat giebt unmittelbar die idealen Primfactoren
der complexen Zahl $\psi_k(\alpha)$.\footnote{This result by Jacobi
immediately gives the ideal prime factors of the complex number 
$\psi_k(\alpha)$.}
\end{quote}
Here Kummer's $\psi_k(\alpha)$ is our $-J(\chi,\chi^k)$, and the result 
by Jacobi alluded to is his fundamental congruence (\ref{FC}); for its 
statement, Kummer refers to the publications \cite{JacKrT} from 1837 and 
1846, and then claims (see \cite[p. 361]{KumZ}) that Jacobi proved this 
``at the place cited''. There are, however, no proofs in \cite{JacKrT}, 
just the statement
\begin{quote}
Die Beweise dieser S\"atze konnten in den vergangenen Wintervorlesungen 
ohne Schwierigkeit meinen Zuh\"orern mitgetheilt werden\footnote{The
proofs of these theorems could be communicated without problems to my 
students in the lectures in the last winter semester.}
\end{quote} 
followed by the footnote directed at Eisenstein that we quoted in 
the Introduction. Since Kummer apparently knew Jacobi's proofs, he 
must have read them in a copy of Jacobi's lectures in number theory.

\section{Kummer's Ideal Numbers}

Let $\lambda$ be a prime and $\alpha$ a primitive root of the equation
$\alpha^\lambda = 1$. The elements of the ring $\Z[\alpha]$ can be
written as polynomials 
$$ f(\alpha) = a_0 + a_1\alpha + \ldots + a_{\lambda-1} \alpha^{\lambda-1} $$
with coefficients $a_j \in \Z$. For an integer $k$ coprime to $\lambda$,
the polynomial $f(\alpha^k)$ is then the conjugate of $f(\alpha)$ 
with respect to the automorphism $\sigma_k: \alpha \too \alpha^k$.
The norm of $f(\alpha)$ is the product 
$f(\alpha) f(\alpha^2) \cdots f(\alpha^{\lambda-1})$.

In \cite{Kum}, Kummer then explains that there are several possible ways
of introducing ideal complex numbers; the simplest idea, and
apparently the one that Kummer came up with first, is based on
the following observation, which Kummer communicated to Kronecker
in a letter\footnote{Kummer's letters to Kronecker can be found
in Kummer's Collected Papers \cite{KCP}.} from April 10, 1844:
\begin{quote}
Wenn $f(\alpha)$ die Norm $p$ hat ($p$ Primzahl $\lambda n + 1$),
so ist jede complexe Zahl einer reellen congruent f\"ur den Modul
$f(\alpha)$. Hierbei ist nur zu zeigen, da\ss{} 
$\alpha \equiv \xi \bmod f(\alpha)$, wo $\xi$ reell. Die\ss{} scheint
sich von selbst zu verstehen, weil $\xi-\alpha$ wenn
$$ 1+ \xi + \xi^2 + \ldots + \xi^{\lambda-1} \equiv 0 \bmod p $$
stets einen Factor mit $p$ gemein hat, wie in dem Beweise, da\ss{} jede
Primzahl $p$ sich in $\lambda-1$ complexe Factoren zerlegen l\"a\ss{}t
gezeigt wird.
\end{quote} 

At this point, Kummer apparently still believed that he had proved that
every prime $p \equiv 1 \bmod \lambda$, where $\lambda$ is an odd prime,
splits into $\lambda-1$ factors in $\Z[\alpha]$, where $\alpha^\lambda = 1$. 
A little later, Jacobi, upon his return from his journey to Italy on
June 17, would point out this mistake to Kummer. Kummer then sat down 
to see what parts of his work survived, and was led to his first attempt 
at defining ideal numbers using Jacobi's maps; in \cite{Kum}, he motivates
his definition by {\em assuming} that the prime $p = \lambda n +1$ 
splits into prime factors in $\Z[\alpha]$: 
$$ p = f(\alpha) f(\alpha^2) f(\alpha^3) \cdots f(\alpha^{\lambda-1}). $$
Kummer then writes
\begin{quote}
Ist $f(\alpha)$ eine wirkliche complexe Zahl und ein Primfactor 
von $p$, so hat sie die Eigenschaft, da\ss, wenn statt der Wurzel 
der Glei\-chung $\alpha^\lambda = 1$ eine bestimmte Congruenzwurzel 
von $\xi^\lambda \equiv 1 \bmod p$ substituirt wird, 
$f(\xi) \equiv 0 \bmod p$ wird.\footnote{If $f(\alpha)$ is an actual 
complex number and a prime factor of $p$, then it has the property 
that when a root of the congruence $\xi^\lambda \equiv 1 \bmod p$ is 
substituted for the root of the equation $\alpha^\lambda = 1$, we get 
$f(\xi) \equiv 0 \bmod p$.}
\end{quote}
The main difference to what he wrote in his letter to Kronecker
in April 1844  is that Kummer now explicitly assumes the
existence of a decomposition of $p$. He then continues:
\begin{quote}
Also auch, wenn in einer complexen
Zahl $\Phi(\alpha)$ der Primfactor $f(\alpha)$ enthalten ist, wird
$\Phi(\xi) \equiv 0 \bmod p$; und umgekehrt: wenn 
$\Phi(\xi) \equiv 0 \bmod p$, und $p$ in $\lambda-1$ complexe
Primfactoren zerlegbar ist, enth\"alt $\Phi(\alpha)$ den Primfactor
$f(\alpha)$. Die Eigenschaft $\Phi(\xi) \equiv 0 \bmod p$ ist nun 
eine solche, welche f\"ur sich selbst von der Zerlegbarkeit der 
Zahl $p$ in $\lambda-1$ Primfactoren gar nicht abh\"angt; sie 
kann demnach als Definition benutzt werden, indem bestimmt wird, 
da\ss{} die complexe Zahl $\Phi(\alpha)$ den idealen Primfactor 
von $p$ enth\"alt, welcher zu $\alpha = \xi$ geh\"ort, wenn 
$\Phi(\xi) \equiv 0 \bmod p$ ist. Jeder der $\lambda-1$ complexen
Primfactoren von $p$ wird so durch eine Congruenzbedingung 
ersetzt.\footnote{Thus if the prime $f(\alpha)$ divides a complex 
number $\Phi(\alpha)$, then we will have $\Phi(\xi) \equiv 0 \bmod p$; 
and conversely: if $\Phi(\xi) \equiv 0 \bmod p$, and $p$ can 
be decomposed into $\lambda-1$ complex prime factors, then 
$\Phi(\alpha)$ contains  the prime factor $f(\alpha)$. The 
property $\Phi(\xi) \equiv 0 \bmod p$ is such that it does not 
depend on the possibility of decomposing the number $p$ into 
$\lambda-1$ prime factors; thus we can use it as a definition, 
by demanding that the complex number $\Phi(\alpha)$ contain the 
ideal prime factor of $p$ belonging to 
$\alpha = \xi$ if $\Phi(\xi) \equiv 0 \bmod p$. Each of the  
$\lambda-1$ complex prime factors of $p$ is replaced by a 
congruence condition in this way.}
\end{quote}
The Jacobi map $\phi:\alpha \too \xi \bmod p$ has the property
that $\phi(f(\alpha)) = f(\xi) \equiv 0 \bmod p$. Thus if
$f(\alpha) \mid \Phi(\alpha)$ in $\Z[\alpha]$, then 
$\Phi(\alpha) = f(\alpha) g(\alpha)$, and applying $\phi$ shows
that 
$\Phi(\xi) = \phi(\Phi(\alpha)) = \phi(f(\alpha)) \phi(g(\alpha))
           = f(\xi) g(\xi) \equiv 0 \bmod p$. Kummer then makes the 
crucial observation that the congruence $\Phi(\xi) \equiv 0 \bmod p$ 
makes sense whether $f(\alpha)$ exists or not: it is a consequence
of the existence of the Jacobi map! Kummer then attaches an ideal 
prime to every Jacobi map $\alpha \too \xi \bmod p$. Then an integer 
$\Phi(\alpha) \in \Z[\alpha]$ will be divisible by the ideal prime
attached to $\phi$ if and only if $\Phi(\xi) \equiv 0 \bmod p$.

Nowadays we would not hesitate {\em defining} an ideal prime 
to be the Jacobi map, but such an idea would probably have been 
too revolutionary even for Dedekind, and certainly must have been 
out of reach for Kummer, who talked about the ideal prime 
{\em belonging to} $\alpha = \xi$ instead. As Kummer explains, 
however, there are problems connected with this approach:
\begin{quote}
In der hier gegebenen Weise aber gebrauchen wir die Congruenzbedingungen
nicht als Definitionen der idealen Primfactoren, weil diese nicht
hinreichend sein w\"urden, mehrere gleiche, in einer complexen
Zahl vorkommende ideale Primfactoren vorzustellen, und weil sie,
zu beschr\"ankt, nur ideale Primfactoren der realen Prim\-zah\-len
von der Form $m\lambda +1$ geben w\"urden.\footnote{We do not use the
congruence conditions in the way given here as definitions of ideal
prime factors since these would not suffice to detect several equal
ideal prime factors occurring in a complex number, and since they,
being too narrow, would only yield ideal prime factors of the real
prime numbers of the form $m\lambda +1$.}
\end{quote}
Thus the problems Kummer was facing were
\begin{enumerate}
\item[(A)] Inertia: the Jacobi maps $\Z[\alpha] \lra \F_p$ only 
      provide ideal numbers dividing primes $p \equiv 1 \bmod \lambda$.
\item[(B)] Multiplicity: there is no obvious way of defining 
      the exact power of an ideal prime dividing a given element
      in $\Z[\alpha]$. 
\item[(C)] Completeness: how can we be sure that we have found ``all''
      ideal primes?
\end{enumerate}
Kummer's solution of these problems will be discussed in the next 
few sections. Afterwards we will explain the close connection 
between Kummer's ideas and modern valuation theory.

It follows immediately from Kummer's definition that ideal primes 
behave like primes: the ideal prime attached to the Jacobi map 
$\phi$ divides $f(\alpha) \in \Z[\alpha]$ if and only if 
$\phi(f(\alpha)) = 0$; if it divides a product $f(\alpha) g(\alpha)$, 
then $0 = \phi(fg) = \phi(f) \phi(g)$, hence it divides a factor.

Before we start addressing the problems listed above, we remark
that the ideal prime dividing $p = \lambda$ is easy to deal with:
there is only one, it is ``real'' (namely $\pi = 1 - \alpha$),
and the corresponding Jacobi map is defined by 
$\phi(\alpha) = 1 + \lambda \Z$.

\section{Solving Problem (A): Decomposition Fields}

The first problem is easy to solve for us: instead of looking at
homomorphisms $\phi: \Z[\alpha] \lra \F_p$, we consider surjective
homomorphisms $\phi: \Z[\alpha] \lra \F_q$ for finite fields with
$q = p^f$ elements. Let $\Phi(X) = 1 + X + X^2 + \ldots + X^{\lambda-1}$
denote the $\lambda$-th cyclotomic polynomial; if 
$\Phi(X) \equiv P_1(X) \cdots P_f(X) \bmod p$ splits into $f$ irreducible 
factors $P_j(X)$ over $\F_p$, then reduction modulo $p$ gives us 
surjective ring homomorphisms
$\phi_j: \Z[\alpha] \simeq \Z[X]/(\Phi) \lra \F_p[X]/(P_j) \simeq \F_q$.
These Jacobi maps can then be used to define ideal numbers for general
primes $p$. 

This was, however, not an option for Kummer: although Gauss had already 
introduced residue class fields in the ring of Gaussian integers, some 
of which do have $p^2$ elements, the general theory of finite fields 
originated in the work of Galois, which was completely unknown in Germany
at the time Kummer started working on these problems\footnote{Galois'
work was published in 1846 by Liouville.}. In a development
independent of the work of Galois, Sch\"onemann \cite{Sch} started
studying ``higher congruences'', as the theory of polynomial rings over
$\F_p$ was called at the time, at about the same time Kummer 
worked on ideal numbers.

So how did Kummer proceed then? With hindsight, Kummer's solution 
is simple and ingenious, and in order to explain why it works, we
will use the language of Dedekind's ideal theory. Let $p \ne \lambda$ 
be a prime with order $f$ in $(\Z/\lambda\Z)^\times$; then $p$ splits 
into $e = \frac{\lambda-1}f$ distinct primes of degree $1$ in the 
decomposition field $F$ of $p$ (since $K/\Q$ is cyclic, this is 
the unique subfield of degree $e$ over $\Q$), say 
$p \cO_F = \fp_1 \cdots \fp_e$. These prime ideals $\fp_j$ remain 
inert in $K/F$. Since $\cO_F/\fp_j \simeq \Z/p\Z$, every element of 
$\cO_F$ is congruent to an integer modulo $\fp_j$, hence reduction
modulo $\fp_j$ defines a Jacobi map $\phi_j: \cO_F \lra \Z/p\Z$. 
The Gaussian periods $\eta_1$, \ldots, $\eta_e$ form an integral 
basis of $\cO_F$, so we have $\cO_F = \Z[\eta_1, \ldots, \eta_e]$, 
and in particular there are integers $u_1$, \ldots, $u_e$ such that 
$\eta_1 \equiv u_1$, \ldots, $\eta_e \equiv u_e \bmod \fp_j$. The 
integers $u_1$, \ldots, $u_e$ completely characterize the Jacobi 
map $\phi_j$, and, therefore, the prime ideal $\fp_j$. Thus we 
can avoid the introduction of finite fields at the cost of replacing 
$\cO_K = \Z[\alpha]$ by $\cO_F$.

None of the facts used above were known to Kummer, who proved the 
existence of these integers $u_j$ ab ovo (they have the property 
that the substitution $\eta_i \too u_i$ turns equations into 
congruences modulo $p$), and then could describe the associated 
Jacobi maps using systems of congruences in $\Z$.

For each prime $p$, Kummer attaches an ideal prime to each Jacobi 
map $\cO_F \lra \Z/p\Z$. This does not really solve Kummer's first 
problem: for deciding whether an element $f(\alpha) \in \Z[\alpha]$ 
is divisible by an ideal prime attached to $\phi$, we would like to 
evaluate $\phi(f)$ and thus face the problem that $\phi$ is defined 
on the subring $\cO_F$ of $\Z[\alpha]$, but not on $\Z[\alpha]$. 
Kummer's solution of problem (B), namely defining a factorization 
of cyclotomic integers into powers of ideal prime numbers, also 
provided him with a clever way around having to extend $\phi$ from 
$\cO_F$ to $\Z[\alpha]$. 

Before we turn to Kummer's solution of problem (B), let us address 
question (C): Suppose we have what we believe to be a complete set 
of Jacobi maps from our domains $\Z[\alpha]$ to certain finite 
fields; how we can be sure to have found ``all'' of them? Kummer's 
answer was as follows: he proved the ``fundamental theorem'' that 
$f(\alpha) \mid g(\alpha)$ if and only if each ideal prime divides 
$g(\alpha)$ with at least the same multiplicity with which it divides 
$f(\alpha)$. This property can be formulated in a slightly different
way: Each Jacobi map $\phi$ is defined at $\frac{g(\alpha)}{f(\alpha)}$
if and only if $f(\alpha) \mid g(\alpha)$. Clearly $\phi$ is defined
at $\frac{g(\alpha)}{f(\alpha)}$ if $f(\alpha) \mid g(\alpha)$, since 
then the quotient is an element of $\Z[\alpha]$. The essential criterium
for completeness therefore is the following:
\begin{quote}
If $h \in \Q(\alpha)$ is an element at which every Jacobi map 
$\phi: \Z[\alpha] \lra \F_q$ is defined, then $h \in \Z[\alpha]$.
\end{quote}
As we will see, this is not just a statement on the completeness
of the Jacobi maps, but also on the correct choice of the ring
of integers we are working with -- in our case $\Z[\alpha]$ and
not some smaller ring.

\section{Solving Problem (B): Valuations}
Now let us look at problem (B): defining multiplicity. We will 
immediately discuss the general case of primes $q$ with 
$q^f \equiv 1 \bmod \lambda$. If we think of ideal numbers not, 
as Kummer did, as systems of congruences but as being attached 
to Jacobi maps, then it is not difficult to define when two numbers 
$f(\alpha), g(\alpha) \in \Z[\alpha]$ are divisible by the same 
power of an ideal prime attached to $\phi$: in such a case we 
would expect that this ideal prime can be cancelled in the 
fraction $\frac{f(\alpha)}{g(\alpha)}$, i.e., that there exist 
elements $f'(\alpha), g'(\alpha) \in \Z[\alpha]$ such that 
$\frac fg = \frac{f'}{g'}$ with $\phi(f') \phi(g') \ne 0$.

This idea can be extended immediately: if we can write 
$\frac fg = \frac{f'}{g'}$ with $\phi(g') \ne 0$, we say that 
the ideal number attached to $\phi$ divides $f$ at least as often
as $g$. If we let $v_\phi(f)$ denote the hypothetical exponent
with which the ideal prime attached to $\phi$ divides a number
$f(\alpha)$, then the above definitions will tell us when
$v_\phi(f) < v_\phi(g)$, $v_\phi(f) = v_\phi(g)$, or 
$v_\phi(f) > v_\phi(g)$. The fundamental problem now is
to show that, given elements $f, g \in \Z[\alpha]$, we always
are in exactly one of these three situations: 
\begin{enumerate}
\item[(B')] Given nonzero elements $f(\alpha), g(\alpha) \in \Z[\alpha]$, 
        there exist $f'(\alpha), g'(\alpha) \in \Z[\alpha]$ with 
        $\frac fg = \frac{f'}{g'}$ and $\phi(f') \ne 0$ or $\phi(g') \ne 0$.
\end{enumerate}
Let us say that the Jacobi map $\phi$ is defined at a nonzero 
element $h(\alpha) \in \Q(\alpha)$ if we can write $h = \frac fg$ 
with $\phi(g) \ne 0$. Then (B') can be formulated in the following way:
\begin{enumerate}
\item[(B'')] Given an element $h(\alpha) \in \Q(\alpha)$, a Jacobi map
        $\phi$ is defined at $h$ or at $\frac1h$.
\end{enumerate}

Subrings $R$ of a field $K$ with the property that for every
$h \in K$ we have $h \in R$ or $\frac1h \in R$ are called valuation 
rings; each Jacobi map satisfying property (B'') defines a valuation 
ring in $\Q(\alpha)$.

For stating (B') and (B'') we have assumed the existence of
the exponent $v_\phi(f)$. In order to guarantee its existence
we have to assume that the valuation ring defined by $\phi$ 
has additional properties. Everything we need will follow from 
\begin{enumerate}
\item[(B''')] For each Jacobi map $\phi$ there exists a 
        $\psi(\alpha) \in \Z[\alpha]$ with the following properties:
    \begin{enumerate}
         \item $\phi(\psi(\alpha)) = 0$;
         \item $\phi$ is defined at $\frac{f(\alpha)}{\psi(\alpha)}$ 
               for all $f(\alpha) \in \Z[\alpha]$ with 
               $\phi(f(\alpha)) = 0$;
         \item if $\phi$ is defined at $f(\alpha)/\psi(\alpha)^n$
               for all $n \in \N$, then $f(\alpha) = 0$.
     \end{enumerate}
\end{enumerate}
Such an element $\psi$ is called a uniformizer for (the valuation
defined by) $\phi$. Let us first show that for any $f \in K^\times$
there is an integer $n \in \Z$ such that $\phi$ is defined at 
$f/\psi^n$: write $f = \frac gh$ for $g, h \in R$; if $\phi(h) \ne 0$, 
then we can choose $n = 0$.  If $\phi(h) = 0$, let $n \ge 0$ be the
maximal integer such that $\phi$ is defined at $h/\psi^n$. Then 
$\phi(h/\psi^n) \ne 0$ by the maximality of $n$ and (b), hence
$\phi$ is also defined at $\psi^n/h$. But then $\phi$ is defined
at $\frac gh \psi^n = f/\psi^{-n}$. Using this argument it is
easy to show the existence of $v_\phi(f)$ as well as property (B'').

The properties (B'), (B'') and (B''') will be discussed, within the 
context of the theory of valuations, in the sections below; in 
the rest of this section, let us see how Kummer solved problem (B):
\begin{quote}
Die von mir gew\"ahlte Definition der idealen complexen Primfactoren,
welche im Wesentlichen zwar mit den hier angedeuteten \"ubereinstimmt,
aber einfacher und allgemeiner ist, beruht darauf, da\ss{} sich,
wie ich besonders beweise, immer eine aus Perioden ge\-bil\-de\-te Zahl
$\psi(\eta)$ finden l\"a\ss{}t von der Art, da\ss{}
$$ \psi(\eta) \psi(\eta_1) \psi(\eta_2) \cdots \psi(\eta_{e-1}) $$
(welches eine ganze Zahl ist) durch $q$ theilbar sei, aber nicht
durch $q^2$. Diese complexe Zahl $\psi(\eta)$ hat alsdann immer
die obige Eigenschaft, da\ss{} sie congruent Null wird, modulo $q$,
wenn statt der Perioden die entsprechenden Congruenzwurzeln gesetzt
werden, also $\psi(\eta) \equiv 0 \bmod q$, f\"ur $\eta = u$,
$\eta_1 = u_1$, $\eta_2 = u_2$ etc. Ich setze nun 
$\psi(\eta_1) \psi(\eta_2) \cdots \psi(\eta_{e-1}) = \Psi(\eta)$
und definire die idealen Primzahlen folgenderma\ss{}en:

``Wenn $f(\alpha)$ die Eigenschaft hat, da\ss{} das Produkt 
$f(\alpha) \cdot \Psi(\eta_r)$ durch $q$ theilbar ist, so 
soll dies so ausgedr\"uckt werden: Es enth\"alt $f(\alpha)$ 
den idealen Primfactor von $q$, welcher zu $u = \eta_r$ 
geh\"ort. Ferner, wenn $f(\alpha)$ die Eigenschaft hat, da\ss{} 
$f(\alpha) (\Psi(\eta_r))^\mu$ durch $q^\mu$ theilbar ist, aber
$f(\alpha) (\Psi(\eta_r))^{\mu+1}$ nicht theilbar durch $q^{\mu+1}$,
so soll dies hei\ss{}en:  Es enth\"alt $f(\alpha)$ den zu $u = \eta_r$ 
geh\"origen idealen Primfactor von $q$ genau $\mu$ 
mal.''\footnote{The definition of the ideal complex prime factors
I have chosen, which coincides essentially with the ones sketched
above, but is simpler and more general, is based on the fact that,
as I will prove, there always exists a number $\psi(\eta)$, formed
out of periods, with the property that 
$$ \psi(\eta) \psi(\eta_1) \psi(\eta_2) \cdots \psi(\eta_{e-1}), $$
which is an integer, is divisible by $q$, but not by $q^2$. This
complex number $\psi(\eta)$ has the property above of becoming
congruent to $0 \bmod q$ if we replace the periods by the 
corresponding roots of the congruence, i.e., 
$\psi(\eta) \equiv 0 \bmod q$ for $\eta = u$, $\eta_1 = u_1$, 
$\eta_2 = u_2$ etc. Now I set 
$\psi(\eta_1) \psi(\eta_2) \cdots \psi(\eta_{e-1}) = \Psi(\eta)$
and define the ideal prime factors as follows:
``If $f(\alpha)$ has the property that the product 
$f(\alpha) \cdot \Psi(\eta_r)$ is divisible by $q$, then we shall
express this by saying that $f(\alpha)$ contains the ideal prime
factor of $q$ belonging to $u = \eta_r$. Moreover, if  $f(\alpha)$
has the property that $f(\alpha) (\Psi(\eta_r))^\mu$ is divisible 
by $q^\mu$ without $f(\alpha) (\Psi(\eta_r))^{\mu+1}$ being divisible
by $q^{\mu+1}$, then this shall mean: $f(\alpha)$ contains the
prime ideal factor of $q$ belonging to $u = \eta_r$ exactly $\mu$ 
times.}
\end{quote}

Kummer's element $\psi(\eta)$ is a uniformizer for the ideal prime
attached to the generalized ``Jacobi map'' $\phi: \eta \too u$, 
$\eta_1 \too u_1$, \ldots, $\eta_{e-1} \too u_{e-1}$; its norm 
$N_{F/\Q} \psi(\eta) = \psi(\eta) \Psi(\eta)$ is divisible by the 
prime $q$, but not by $q^2$. The multiplicity $\mu$ with which the
ideal prime attached to $\phi$ divides some $f(\alpha) \in \Z[\alpha]$
is the maximal natural number $\mu$ with the property that 
$f(\alpha) \Psi(\eta)^\mu$ is divisible by $q^\mu$. Rephrasing
this condition slightly, we see that it is equivalent to the
condition that $\phi$ is defined at 
$\frac{f(\alpha) \Psi(\eta)^\mu}{q^\mu}$, or, equivalently
(observe that $\phi$ is defined at both $\frac{\psi(\eta)}q$ 
and its inverse), at $\frac{f(\alpha)}{\psi(\eta)^\mu}$.

\section{Jacobi Maps and Valuations}
Kummer developed a ``complete'' theory of ideal numbers only for
cyclotomic fields and their subfields. Although he ran into
insurmountable problems when trying to extend his construction
to general number fields, he did not investigate exactly which 
properties of the rings he was working in were responsible for
his success (it is no exaggeration to claim that no one before
Dedekind understood such questions properly). 

In order to display the gaps in Kummer's theory as clearly
as possible we now try to transfer his construction to more
general situations. Let $R$ be a domain with quotient field $K$.
A Jacobi map is a ring homomorphism $\phi: R \lra F$ onto a 
field $F$. The kernel $\fp = \ker \phi$ satisfies $R/\fp \simeq F$,
hence is a maximal ideal (and a fortiori a prime ideal) in $R$;
Jacobi maps with the same kernel will be identified.

A Jacobi map $\phi$ is said to be defined at $t \in K$ if there
exist $r, s \in R$ with $t = \frac rs$ and $\phi(s) \ne 0$. The
set of all $t \in K$ at which $\phi$ is defined is a subring of $K$
denoted by $R_\phi$. 

A subring $R$ of a field $K$ is called a valuation ring if for any 
$t \in K$ we have $t \in R$ or $\frac1t \in R$. The unit group of
a valuation ring consists of all elements $a \in K^\times$ for which 
$a \in R$ and $\frac1a \in R$. The set $\fm = R \setminus R^\times$ 
of nonunits is an ideal, and in fact the unique maximal ideal of 
$R$.\footnote{For the history of valuation theory, beginning in 1912 
with K\"urschak's work, see Roquette's article \cite{Roq}.}

If $\phi$ is a Jacobi map defined on $R$, then $R_\phi$ is a 
valuation ring if and only if the analog of property (B'') holds,
that is: given any $t \in K$, a Jacobi map $\phi$ is 
defined at $t$ or at $\frac1t$. 

An additive valuation\footnote{If $v$ is an additive valuation
and $c$ any real number $> 1$, then $|\alpha| = c^{-v(\alpha)}$ 
defines a ``multiplicative valuation'', or an ``absolute value''
on $K$ (we put $|0| = 0$).} is a map from the nonzero elements of 
a field $K$ to an ordered group $G$ with the properties
\begin{enumerate}
\item $v(ab) = v(a) + v(b)$, 
\item $v(a+b) \ge \min \{v(a),v(b)\}$ 
\end{enumerate}
for all $a, b \in K$. We also set $v(0) = \infty$ and 
$\infty \ge g$ for all $g \in G$; then $v$ is a map 
$K \lra G \cup \{\infty\}$. The set of all $a \in K$ with 
$v(a) \ge 0$ forms a valuation ring. For example, let $p$ be a 
prime and let $v_p(a)$ denote  the exponent of $p$ in the prime 
factorization of the nonzero rational number $a \in \Q$; then 
$v_p(a)$ is a valuation.

Thus every valuation determines a valuation ring (valuations giving rise
to the same valuation ring are called equivalent). Conversely, to every
valuation ring we can find a corresponding valuation: to this end, 
define an order on the ``divisibility group'' $G_K = K^\times/R^\times$
(in \cite{LEuk}, this group shows up in a different connection)  
via $y R^\times \le x R^\times$ if and only if $\frac xy \in R$. It is 
a simple exercise to show that the map $v: K^\times \lra G_K$ sending 
$x \in K^\times$ to its coset $x R^\times$ is indeed a valuation.

For solving Kummer's problem (B) of defining multiplicity, a valuation 
is not good enough: we want a valuation with values not in some ordered 
group but in $\Z$! Such valuations are called discrete valuations; they 
can be characterized by the fact that the corresponding valuation rings 
are discrete valuation rings, i.e., their maximal ideal $\fm$ must be 
principal. In fact, if $\fm = (t)$ is generated by an element 
$t \in K^\times$ (such elements are called uniformizers for the 
corresponding valuation), then any $a \in K^\times$ can be written 
uniquely in the form $a = ut^m$ for some  unit $u \in R^\times$ and 
some $m \in \Z$, and setting $v(a) = m$ gives us a discrete valuation.

In Kummer's case of cyclotomic rings of integers, the existence
of a uniformizing element $h(\alpha)$ is exactly what we asked 
for in property (B'''). Thus in terms of modern algebra, Kummer 
used Jacobi maps to construct discrete valuation rings.

\section{Examples of Jacobi Maps}

The major sources of Jacobi maps and their associated valuations 
are number theory and algebraic geometry. In addition to the two 
most important examples coming from number fields and algebraic
curves, we will also discuss certain monoids that behave a lot
like them.

\subsection*{Hilbert Monoids}

In his lectures on number theory in the winter semester 1897/98,
Hilbert \cite{Hil} used the monoid $\{1, 6, 11, 16, 21, \ldots\}$ 
of natural numbers $\equiv 1 \bmod 5$ for motivating the 
introduction of ideals. We would like to show now that these
monoids can also be used for explaining Kummer's solution of 
problem (B), namely defining multiplicity. A general 
theory of divisibility in monoids was given by Rychlik \cite{Ry}.

Consider the monoid $M = \{1, 5, 9, 13, \ldots\}$ of natural 
numbers $\equiv 1 \bmod 4$. In this monoid, factorization into 
irreducibles is not unique, as the example $21 \cdot 21 = 9 \cdot 49$ 
shows. One can restore unique factorization by adjoining ``ideal 
numbers'' such as $\gcd(9,21)$ representing a common divisor of $9$ 
and $21$. 

Kummer's approach using Jacobi maps also works here: For each
odd prime $p \in \N$, there is a surjective homomorphism of
monoids $\phi_p: M \lra \F_p$. For primes $p \equiv 1 \bmod 4$,
the kernels of the $\phi_p$ are ``principal'' in the sense that
$\ker \phi_p = pM$. For primes $p \equiv 3 \bmod 4$, this is 
not the case: $\ker \phi_3 = \{9, 21, 33, 45, \ldots\}$ cannot
be written in the form $aM$ for some $a \in M$, but we can think
of the Jacobi map $\phi_3$ as representing an ``ideal prime'' $3$ 
in $M$.

In order to define multiplicity we extend the Jacobi maps 
$\phi_p$ to the quotient monoid $Q(M)$ of all fractions $\frac ab$ 
with $a, b \in M$ in the following way: if we can write 
$\frac ab = \frac cd$ with $c, d \in M$ and $\phi_p(d) \ne 0$,
then we set $\phi_p(\frac ab) = \frac{\phi_p(c)}{\phi_p(d)}$ and say
that $\phi_p$ is defined at $\frac ab$; we also set 
$\phi_p(\frac ab) = \infty$ if  $\phi_p(\frac ba) = 0$. 
Since $\frac{9}{21} = \frac{9\cdot 49}{21 \cdot 49} = 
                      \frac{21 \cdot 21}{21 \cdot 49} = \frac{21}{49}$,
we see that $\phi_3$ is defined at $\frac{9}{21}$, and that 
$\phi_3(\frac{9}{21}) = 0$. Similarly, $\phi_3$ is defined at
$\frac{9}{21^2} = \frac{9}{9 \cdot 49} = \frac1{49}$.

The following fact is a fundamental property of $M$: if $\phi_p$
is not defined at $\frac ab$, then it is defined at $\frac ba$.
In fact, if both $c$ and $d$ are divisible in $\N$ by $p$, then
we can cancel $p$ immediately if $p \equiv 1 \bmod 4$, or we 
can find a natural number $r \equiv 3 \bmod 4$ coprime to $p$
and write $\frac cd = \frac{pc'}{pd'} = \frac{rc'}{rd'}$ with 
$rc', rd' \in M$. Continuing in this way we see that $\frac ab = \frac cd$
with $p \nmid c$ or $p \nmid d$, and then $\phi_p$ is defined at
$\frac ba$ and $\frac ab$, respectively.  

An element $q \in M$ is called a uniformizer for $p$ if $\phi_p$
is defined at $\frac cq$ for all $c \in M$ with $\phi_p(c) = 0$.
It is easy to see that the elements $21, 33, \ldots$ are uniformizers 
for $p = 3$, and that $9$ is not. More generally, given a prime 
$p \equiv 3 \bmod 4$, any element $q = pr$, where $r \equiv 3 \bmod 4$ 
is coprime to $p$, is a uniformizer for $p$.

Now we say that the ideal prime $p$ divides some $a \in M$ with
exponent $m$ if and only if $\phi_p$ is defined at $\frac{a}{q^m}$,
but not at $\frac{a}{q^{m+1}}$. Such an integer $m$ exists, and 
it does not depend on the choice of the uniformizer $q$. Since e.g. 
$\phi_3$ is defined at $\frac{9}{21^2}$, but not at $\frac{9}{21^3}$, 
we find that $9$ is exactly divisible by the square of the ideal 
prime $3$.

It can be shown that every element in $M$ can be written 
uniquely as a product of powers of ideal primes. We leave it as 
an exercise for the reader to show that unique factorization into 
ideal primes in $M$ can be used to show that $a \in M$ is a square 
in $Q(M)$ if and only if it is a square in $M$.

We have seen above that the ``Hilbert monoid'' $M$ is quite
well behaved in that it satisfies the analog of property (B''): 
for odd primes $p$ and all $a, b \in M$, the ``Jacobi maps'' 
$\phi_p$ are defined either at $\frac ab \in Q(M)$ or at 
$\frac ba$ (or at both). 

In fact, the same thing holds in general monoids of Hilbert type: 
for any natural number $m > 1$, let $H$ be a subgroup of 
$G = (\Z/m\Z)^\times$, and define the monoid $M^G_H$ as the set of 
all natural numbers $a \in \N$ whose residue classes mod $m$ lie 
in $H$. The monoid $M = \{1, 5, 9, \ldots\}$ considered above, for 
example, is $M^G_H$ for the trivial subgroup $H$ of $G = (\Z/4\Z)^\times$.
It is easy to show that all such monoids $M^G_H$ are nonsingular in 
the sense that Jacobi maps $\phi_p$, for all primes $p$ coprime
to $m$, are defined at $\frac ab \in Q(M^G_H)$ or at $\frac ba$. 

It is actually not difficult to prove that all these monoids have 
unique factorization into ideal primes, to define an analog 
$\Cl(M^G_H)$ of the ideal class group\footnote{Call two ideal numbers 
$p, q$ equivalent if there is an ideal number $r$ such that 
$pr, qr \in M^G_H$.}, and to show that\footnote{There is an exact sequence
$$ \begin{CD}
   1 @>>> H @>>> G @>{\pi}>> \Cl(M^G_H) @>>> 1,
   \end{CD} $$
where $\pi$ maps the residue class $a + m\Z$ to the ideal
class generated by the ideal number $a$.}  $\Cl(M^G_H) \simeq G/H$.

\subsection*{Number Fields}
In a Dedekind domain $R$, every prime ideal $\fp \ne (0)$ is maximal,
hence induces a Jacobi map $\phi_\fp: R \lra F = R/\fp$. These maximal
ideals also give rise to valuations $v_\fp$: given any $a \in K^\times$, 
write $(a) = \fp^m \fa$ for some ideal $\fa$ in whose prime ideal 
factorization $\fp$ does not occur, and set $v_\fp(a) = m$. 
The valuations $v_\fp$ attached to nonzero prime ideals are non-archimedean
in the sense that the corresponding absolute values $|\,\cdot\,|_\fp$ 
defined by $|\alpha|_\fp = (N\fp)^{-v_\fp(\alpha)}$ satisfy the 
strong triangle inequality 
$|\alpha + \beta|_\fp \le \max \{|\alpha|_\fp, |\beta|_\fp\}$. 
It can be shown that every archimedean valuation in a number
field $K$ is equivalent to a valuation $v_\fp$ for a suitable
prime ideal $\fp$.  

In addition, every number field has some archimedean valuations:
if $K = \Q(\alpha)$ and $\alpha$ is the root of the monic polynomial
$f \in \Z[X]$ of degree $n$, let $\alpha_1 = \alpha$, $\alpha_2$, 
\ldots, $\alpha_n$ denote the roots of $f$ in $\C$. Every element
of $K$ can be written as a polynomial in $\alpha$, say as $g(\alpha)$
with $g \in \Q[X]$; therefore we can define an absolute value 
$|\cdot|_j$ on $K$ by setting $|g(\alpha)|_j = |g(\alpha_j)|$. It
turns out that complex conjugate roots give rise to the same 
absolute value, and that there are $r+s$ independent absolute 
values if the number of real and complex roots of $f$ is $r$ 
and $2s$, respectively (complex roots come in pairs since $f$ 
has real coefficients).

The main difference between archimedean and non-archimedean
valuations from our point of view is that non-archimedean
valuations come from (resp. give rise to) an additive valuation
and thus to a valuation ring, whereas archimedean valuations do not.

\subsection*{Algebraic Curves}
Let $K$ be an algebraically closed field, and $f \in K[X,Y]$ an
irreducible polynomial. The zero set $\cC(K) = \{(x,y) \in K \times K:
 f(x,y) = 0\}$ of $f$ is called a plane algebraic curve. Its
coordinate ring is defined to be the ring $\cO = K[\cC] = K[X,Y]/(f)$, 
whose elements are polynomials in $x = X + (f)$ and $y = Y + (f)$; 
since $f$ was assumed to be irreducible in the factorial domain 
$K[X,Y]$, the coordinate ring of $\cC$ is actually a domain, and 
its quotient field $K(\cC)$ is called the function field of $\cC$.

For each point\footnote{In order to keep things as simple as
possible, we only consider the affine plane here, and remark 
that, in algebraic geometry, the projective point of view is
usually to be preferred. We also have to admit that our definition
of a plane algebraic curve is quite naive; but it suffices for our 
purposes.} $P \in \cC(K)$, we can define a ``Jacobi map'' 
$\phi_P: \cO \lra K$ via $\phi_P(g) = g(P)$. Since $\cO$ 
contains all the constant functions, the $\phi_P$ are surjective 
ring homomorphisms. An element $\frac gh \in K(\cC)$ is said to 
be defined at a point $P \in K \times K$ if $\frac gh = \frac rs$ 
for $r,s \in \cO$ with $\phi_P(s) = s(P) \ne 0$. 

The functions from $K(\cC)$ defined at a fixed point $P$ form a 
ring $\cO_P$ called the local ring at $P$, and $\cO_P$ clearly
contains $K$ (here we are identifying $K$ with the constant 
functions, which are defined everywhere) and even $\cO$. 

As an example, consider the curve defined by $Y^2 = X^3 + X^2$.  
The function $h(x,y) = \frac{y}{x+1} \in K(\cC)$, where $x = X + (f)$ 
and $y = Y + (f)$, is not defined at $P = (-1,0)$, but its inverse 
$\frac{x+1}y$ is: we have 
$\frac{x+1}y = \frac{(x+1)y}{y^2} = \frac{(x+1)y}{x^3+x^2} = \frac{y}{x^2}$,
and this function $\frac1h$ is defined at $P$ and vanishes there.

\section{Singularities}
After having gone through several examples where Jacobi maps
$\phi: R \lra F$ give rise to valuation rings and, therefore, to 
valuations, we will now present examples of Jacobi maps $\phi$
for which $R_\phi$ is not a valuation ring.

\subsection*{Singular Monoids}
In Hilbert monoids $M^G_H$, every Jacobi map $\phi: M^G_H \lra \F_p$
gives rise to a ``valuation''. The situation is completely different for the 
monoid $N = \{1, 2, 4, 5, 6, 8, 9, \ldots\}$ of natural numbers congruent 
to $0, 1, 2 \bmod 4$. Here, the Jacobi map $\phi_2: N \lra \F_2$ cannot 
be extended to the quotient module $Q(N)$ by imitating the process that 
worked so well for $M^G_H$: the map $\phi_2$ is undefined both for 
$\frac 62$ and $\frac 26$. In fact, if we had e.g.  $\frac62 = \frac ab$ 
for some odd $b \in N$, then $b \equiv 1 \bmod 4$; on the other hand, 
$6b = 2a$ implies $a = 3b$ in $\N$, and thus $a \equiv 3 \bmod 4$, 
contradicting the assumption that $a \in N$. 

The element $\frac62 \in Q(N)$ displays a singular behavior in more
than one way: the element $9 = (\frac 62)^2$ is a square in the quotient 
monoid $Q(N)$, but not in $N$. In the natural numbers $\N$, on the other 
hand, square roots are either integers or irrational. The following proof 
of the irrationality of $\sqrt{m}$ for nonsquares $m$ brings out a 
fundamental concept that we will need to explain Kummer's success at 
introducing ideal numbers in cyclotomic rings of integers:
\begin{align*}
 \sqrt{m} \text{ is rational } 
      & \iff \text{ the roots of } X^2 - m \text{ are rational} \\
      & \iff X^2 - m \text{ factors over } \Q  \\
      & \iff X^2 - m \text{ factors over } \Z 
\end{align*}
Here, the only nontrivial step is the last equivalence, which is 
a special case of

\medskip\noindent
{\bf Gauss's Lemma.} 
{\em Let $f$ be a monic polynomial with integral coefficients. Then $f$
     factors over $\Q[X]$ if and only if it factors over $\Z[X]$.}
\medskip

As a matter of fact, Gauss's Lemma (which can be found in Gauss's
Disquisitiones \cite{Gau}) implies the following general result:
algebraic integers (roots of monic polynomials $\in \Z[X]$) are either 
elements of $\Z$ or irrational.

Although it means stretching the analogy beyond its limit, let us 
now talk about ``polynomials'' with coefficients in the monoid $N$
(we can do so since $N$ is a subset of $\N$; note, however, that $N$
is not closed with respect to addition). Then the natural analog of 
Gauss's Lemma does not hold in $N$: the polynomial $f(X) = X^2-9$ is 
irreducible over $N$ in the sense that $f$ cannot be written as a 
product of two linear factors with coefficients in $N$, whereas 
$X^2 - 9 = (X - \frac62)(X + \frac62)$ factors over the quotient 
monoid $Q(N)$.

\subsection*{Singular Curves}
Consider the plane algebraic curve defined by the polynomial 
$f(X,Y) = Y^2 - X^3 - X^2$. Then  both $\frac xy \in K(\cC)$ and its
inverse $\frac yx$ are undefined at the singular point $P = (0,0)$. 
In fact, assume that $\frac xy = \frac gh$ for $g,h \in \cO$ with 
$h(P) \ne 0$. Then $xh(x,y) = yg(x,y)$, that is, 
$X h(X,Y) - Yg(X,Y) = r(X,Y)f(X,Y)$ for polynomials $g, h, r \in K[X,Y]$.
Plugging in $Y = 0$ yields $Xh(X,0) = -(X^3+X^2)r(X,0)$, and cancelling
$X$ gives $h(X,0) = -(X^2+X)r(X,0)$. Plugging in $X = 0$ now shows that 
$h(P) = 0$. A similar calculation shows that $\frac yx$ is not defined at $P$.

As in the case of the monoid $N$, the singular behavior of $P = (0,0)$
is connected with the failure of Gauss's Lemma: the polynomial
$F(T) = T^2 - x-1 \in R[T]$ is irreducible over $R$, but factors as
$F(T) = (T - \frac yx)(T + \frac yx)$ over the field of fractions
$K = k(x,y)$ of $R$. Similarly, the element $x+1 \in \cO = K[\cC]$ 
is not a square in $\cO$, but becomes a square in the quotient field
$K(\cC)$ of $\cO$.

The notion of singularity has its origin in the theory of algebraic
curves, and is connected with the existence of tangents. In fact, 
let $P = (a,b)$ be a point on some plane algebraic curve $\cC$ defined
by a polynomial $f \in K[X,Y]$, and let $f_X$ and $f_Y$ denote the partial 
derivatives of $f$ with respect to $X$ and $Y$. Then the tangent to $\cC$ 
at $P$ is defined to be the line $f_X(P)(X-a) + f_Y(P)(Y-b) = 0$; points
with $f_X(P) = f_Y(P) = 0$ are called singular points. 

In our example of the curve $\cC$ defined by 
$f(X,Y) = Y^2 - X^3 - X^2 = 0$, the only common solution of the two 
equations $f_X = -3X^2 - 2X = 0$ and $f_Y = 2Y = 0$ is $P = (0,0)$, 
which therefore is the only singular point on $\cC$ in the affine plane. 

It can be shown that $\cO_P$ is a valuation ring if and only
if $P$ is a nonsingular point of $\cC$. The connection between
the definition of singularity above and the failure of Gauss's
Lemma in the coordinate ring is provided by the crucial observation
that Gauss's Lemma for monic polynomials holds\footnote{We say that 
Gauss's Lemma holds over a domain $R$ with quotient field $K$ if monic 
polynomials in $R[X]$ that factor in $K[X]$ also factor in $R[X]$.} 
in a domain $R$ if and only if $R$ is integrally closed 
(see e.g. \cite{LEuk}).

Recall that a domain $R$ with quotient field $K$ is called integrally 
closed\footnote{This is a notion introduced by Emmy Noether, modeled 
on Dedekind's definition of an algebraic integer.} if every $a \in K$ 
that is a root of a monic polynomial in $R[T]$ actually belongs to $R$. 
The coordinate ring $R = K[\cC]$ in the example above is not integrally 
closed: the root $\frac yx$ of the monic polynomial $T^2-x-1$ is an 
element of $K$ that is integral over $R$, yet does not belong to $R$.

An irreducible algebraic curve is called normal\footnote{This 
definition also applies to algebraic varieties of higher dimension.} 
if its coordinate ring $K[\cC]$ is integrally closed. Normality is 
a local property: $\cC$ is normal if and only if all its local rings 
$\cO_P$ are integrally closed (\cite[Chap. II, \S\ 5]{Sha}). Finally, 
an algebraic curve is normal if and only if none of its points is 
singular.

\subsection*{Singular Orders}
Consider the ring $\cO = \Z[\sqrt{-3}\,]$; the Jacobi map
$\phi_2: \cO \lra \F_2$ defined by 
$\phi_2(a + b\sqrt{-3}\,) \equiv a+b \bmod 2$
is a surjective ring homomorphism, and thus should correspond
to an ideal prime. Since $\phi_2(1+\sqrt{-3}\,) = \phi_2(2) = 0$,
we expect that $\phi_2$ is defined at $\frac{1+\sqrt{-3}}2$ or
$\frac{2}{1+\sqrt{-3}}$, but a simple calculation shows that this
is not the case. This shows that Kummer's idea of attaching an
ideal prime to each Jacobi map does not work here.

Dedekind's ideal theory, by the way, also fails in $\cO$: consider 
the ideal $\ker \phi_2 = \fp = (2,1+\sqrt{-3}\,)$ in $\cO$; since 
$\cO/\fp$ has two elements, $\fp$ is prime (even maximal). Yet
$\fp^2 = (4,2+2\sqrt{-3}, -2+2\sqrt{-3}\,) = (2)(2,1+\sqrt{-3}\,) = (2)\fp$,
and if we had unique factorization into prime ideals, this would imply
$\fp = (2)$, which is not true, however.

Similar examples are provided by the subrings 
$\cO = \Z[pi] = \{a+pbi: a, b \in \Z\}$ of $\Z[i]$, where $p$ is any 
rational prime. Here $\phi(a+pbi) = a + p\Z$ defines a Jacobi map 
$\cO \lra \F_p$ with kernel $\fp = \ker \phi = p\Z \oplus pi\Z$. 
Since $\cO/\fp \simeq \F_p$, $\fp$ is a prime ideal. On the other 
hand, $\phi$ is not defined at $\frac{pi}p$ and $\frac{p}{pi}$.

It should not be surprising that Gauss's Lemma also fails in 
$R = \Z[\sqrt{-3}\,]$: the polynomial $T^2+T+1$ is irreducible 
over $R$, but factors as  $T^2+T+1 = (T - \rho)(T - \rho^2)$ 
over its field of fractions $K = \Q(\sqrt{-3}\,)$, where 
$\rho = \frac{-1+\sqrt{-3}}2$.

Orders $\cO$ are subrings of the rings of integers $\cO_K$ of a 
number field containing $\Z$; a typical example is $\cO = \Z[\sqrt{-3}\,]$. 
The ``bad'' prime ideals (those attached to Jacobi maps not
satisfying condition (B')) are those dividing the conductor of the
order; the conductor is an ideal that measures how far the order 
is from being maximal: the maximal order $\cO_K$, for example, is 
nonsingular and has conductor $(1)$. The order $\Z[\sqrt{-3}\,]$,
on the other hand, has conductor $(2)$, and the prime ideal dividing
$(2)$ shows a singular behavior.

\section{Integral Closure}

We now would like to use a given Jacobi map $\phi: R \lra F$ for defining
a valuation on $R$. As we have seen, there are situations in which 
this does not work, for example if the domain $R$ is not integrally 
closed. It is actually not hard to show that valuation rings $R$ are 
integrally closed: if $t \in K$ is integral over $R$, then 
$t^n = \sum_{j=0}^{n-1} a_j t^j$. If $t \not\in R$, then 
$\frac1t \in R$ since $R$ is a valuation ring; this implies
$t = \sum_{j=0}^{n-1}a_j t^{j+1-n}$. But then $t \in R$, which is
a contradiction.

Assume now that $K$ is a number field with ring of integers $\cO_K$,
and that $\cO \subset \cO_K$ is a proper subring containing $\Z$
(such subrings are called orders). Then there is an integral element 
$\alpha \in K \setminus \cO$. Since $\alpha$ is not in $\cO$, there 
must be a Jacobi map $\phi: \cO \lra \F_q$ such that $\phi$ is not 
defined at $\alpha$. We claim that $\phi$ is also not defined at 
$\frac1\alpha$.

Suppose that we can write $\alpha = \frac{\beta}{\gamma}$ with 
$\phi(\beta) \ne 0$. From $\alpha^n = \sum_{j=0}^{n-1} a_j\alpha^j$ for 
integral coefficients $a_j \in \Z$ we deduce, after multiplying through 
by $\gamma^n$, that 
$\beta^n = \gamma(a_{n-1}\beta^{n-1} + \ldots + a_0 \gamma^{n-1})$.
This implies that $\phi(\gamma) \ne 0$. But then $\phi$ is defined 
at $\alpha = \frac{\beta'}{\gamma'}$ contrary to our assumption.

Thus if $\cO$ is not integrally closed, then there exist elements 
$\alpha \in \cO_K \setminus \cO$ and Jacobi maps $\phi: \cO \lra F$ 
which are not defined at $\alpha$ and $\frac1\alpha$. In these cases, 
Kummer's method of attaching an ideal prime to Jacobi maps fails.

\section{Ideal Numbers and Integral Closure}

Although Kummer did not isolate the property (B') (let alone (B'') 
or even (B''')) in his work, for verifying that his prime ideal 
exponents have the desired properties he implicitly had to 
prove some form of property (B'). In fact, we have seen above that 
Kummer's method of attaching an ideal prime to each Jacobi map 
sometimes does not work. The reason why Kummer did not run into 
these problems, at least not until he had to study ideal numbers 
in Kummer extensions of cyclotomic fields (these are extensions 
of $K = \Q(\zeta_p)$ of the form $L = K(\sqrt[p]{\mu}\,)$; for the 
rings $\cO \subset L$ that Kummer considered, he had to exclude all 
ideal primes dividing $(1-\zeta_p)\mu$), was clarified much later 
by Dedekind. Dedekind was the first to give the correct definition of 
algebraic integers\footnote{Dirichlet, for example, proved his unit 
theorem in orders $\Z[\alpha]$, where $\alpha$ is an algebraic integer, 
i.e., a root of a monic polynomial with integral coefficients.

On the other hand E.~Heine defined algebraic integers in \cite{Heine} 
as numbers that can be constructed from the rational integers by 
addition, multiplication, and raising to $m/n$-th powers with $m, n$ 
positive integers. He then goes on to show that every such number is 
integral in the modern sense, i.e. it is a root of a monic polynomial 
with integral coefficients. The converse is, of course, false, as
Heine's construction gives only algebraic numbers that are
{\em solvable}, i.e., that can be expressed in terms of radicals;
Heine  claimed that any root of a solvable monic polynomial with
integral coefficients is integral in his sense, but the proof he 
gave is not valid.}, and showed that the ring $\cO_K$ of all algebraic 
integers contained in a number field $K$ is ``nonsingular'' in the 
sense above. When Emmy Noether later characterized ``Dedekind domains''
(these are domains in which every ideal can be written uniquely as a
product of prime ideals) axiomatically, integral closure turned out
to be one of the axioms.

Luckily for Kummer, the obvious choice of the ring $\Z[\alpha]$ in 
a cyclotomic field $K = \Q(\alpha)$ turns out to be the full ring
$\cO_K$ of integers. But Kummer found the correct ring of integers
even for the subfields of cyclotomic fields: the rings 
$\Z[\eta_1, \ldots, \eta_e]$ generated by the Gaussian periods 
$\eta_i$, which Gauss had introduced in the seventh section of 
his Disquisitiones Arithmeticae \cite{Gau}! Note that for quadratic
number fields $\Q(\sqrt{p}\,)$ with $p \equiv 1 \bmod 4$, the 
ring of periods is $\Z[\frac{1+\sqrt{p}}2\,]$. 

Kummer mentioned the ring of periods already in his letter to Kronecker 
from October 2, 1844; the main part of his letter was, however, devoted 
to a proof that $\Z[\zeta_5]$ is a Euclidean ring. But although Kummer 
suspected early on that the solution to his problems could be found in 
this ring of periods, it took him a whole year to work out the
details: in his letter from October 18, 1845, he finally could
explain his new theory of ideal numbers to Kronecker.

After having discussed the relevance of integral closure for Kummer's
theory in the last few sections, we are left with the following question: 
where exactly did Kummer use the fact that the rings $\Z[\alpha]$ and 
$\Z[\eta_1, \ldots, \eta_e]$ are integrally closed? Kummer's first 
construction of ideal numbers contained a serious gap (already 
noticed by Eisenstein, and later by Cauchy and Dedekind), which is why 
Edwards presents Kummer's ``second'' proof in \cite{EdwF}; but even 
there it is not at all clear whether Kummer actually used the integral 
closure of the rings he considered. Simply going through his claims with 
the example $\cO = \Z[\sqrt{-3}\,]$ and the ideal prime attached to the
Jacobi map $\phi: \cO \lra \Z/2\Z$ sending $1$ and $\sqrt{-3}$ to the 
residue class $1 + 2\Z$ does not help very much: Kummer's 
construction is so tied to special properties of the rings he
is working in that one has to make various choices when attempting to
transfer his theory to general number fields, and the places where 
integral closure is needed depends on the choices that are made.

Consider, for example, the question whether there exist elements 
$\Psi \in \cO$ such that $\alpha \cdot \Psi \equiv 0 \bmod 2$ if 
and only if $\phi(\alpha) = 0$ for $\alpha \in \cO$. It is easy 
to verify that $\Psi = 1+\sqrt{-3}$ has this property, but this
is not the element you get by faithfully imitating Kummer's 
construction. In particular, our $\Psi$ has the property
$\Psi^2 \equiv 0 \bmod 2$, whereas Kummer proves and uses 
the fact that $\Psi^2 \not\equiv 0 \bmod p$ for his elements. 

In his dissertation \cite{Haub} on the genesis of Dedekind's ideal 
theory, Haubrich also discussed Kummer's construction. A look into 
his very informative thesis quickly reveals that even Kummer's 
second proof contains a gap! This gap was first noticed by Dedekind 
in his review \cite[p. 418]{Ded} of Bachmann's book \cite{BachK}: 
Kummer (as well as Bachmann) did not prove that if an ideal 
prime divides a number $m$ times it also divides it $m-1$ times; 
Dedekind remarked that, as long as this has not been accomplished, 
it is conceivable that an ideal prime divides a number exactly six 
times and exactly eight times. Haubrich also explains that Dedekind's 
proof of the corresponding fact for ideals uses the integral closure 
of the domain he was working in.

\section{Summary}
We have seen that Jacobi introduced techniques which, in our language,
give rise to ring homomorphisms $\phi$ from the ring $\Z[\alpha]$ of 
the ring of $\lambda$-th roots of unity to $\Z/p\Z$ for primes 
$p \equiv 1 \bmod \lambda$. These maps, as Kummer realized, could be 
used for defining ``ideal numbers''. Kummer generalized the Jacobi maps
to all primes $p$ and characterized them by certain sets of congruences
(this is completely in line with the spirit of Kronecker, whose 
ultimate goal was to reduce mathematics to working with natural 
numbers). While these congruences were easy to work with in practice, 
they also made the algebraic structure behind the construction of
ideal numbers invisible, and Dedekind had to struggle for  quite a 
while before he could uncover this structure again: he arrived at 
his ideals by looking at the set of all algebraic integers of a 
number field divisible by Kummer's ideal numbers. Had he started 
with Jacobi maps instead of Kummer's sets of congruences, he might 
have arrived earlier at the correct definition of ideals: these are 
exactly the kernels of Jacobi's ring homomorphisms $\phi$.

Expositions of Kummer's theory of ideal numbers were given
by Bachmann \cite{BachK} and Edwards \cite{EdwF}. Both authors
give Kummer's description of ideal numbers in terms of sets of
congruences; Edwards (\cite[Sect, 4.9, Ex. 10]{EdwF}) gives 
the bijection between ideal numbers and ring homomorphisms 
$\Z[\alpha] \lra \F_q$ as an exercise, and mentioned the
connection to valuation theory in \cite{Edw4}. Dieudonn\'e, in 
his review (MR1160701) of \cite{Edw4}, also recognized the 
algebraic structure behind Kummer's construction when he wrote
\begin{quote}
what [\ldots] Kummer did was the determination of discrete 
valuations on a cyclotomic field.
\end{quote}

Nevertheless, the simple algebraic idea behind Kummer's construction 
remained almost as unknown as Jacobi's role in the creation of ideal
numbers.

\section*{Acknowledgments}
I thank Harold Edwards for his comments on a preliminary draft
of this article.

\end{document}